\documentclass[11pt]{article}
\usepackage[margin=1in]{geometry}
\usepackage{amsmath,amssymb}
\usepackage{mathtools}
\usepackage{tikz}
\usetikzlibrary{arrows.meta}

\mathtoolsset{showonlyrefs}
\usepackage{array,booktabs,longtable}

\usepackage[
  colorlinks=true,
  linkcolor=blue,
  citecolor=blue,
  urlcolor=blue
]{hyperref}
\numberwithin{equation}{section}
\usepackage{amsthm}
\theoremstyle{plain}
\newtheorem{theorem}{Theorem}[section]
\newtheorem{corollary}[theorem]{Corollary}
\newtheorem{proposition}[theorem]{Proposition}
\newtheorem{lemma}[theorem]{Lemma}
\theoremstyle{definition}
\newtheorem{definition}[theorem]{Definition}
\theoremstyle{remark}
\newtheorem{remark}[theorem]{Remark}
\usepackage{todonotes}
\usepackage[maxnames=10,giveninits=true]{biblatex}
\DeclareFieldFormat[article]{title}{\mkbibemph{#1}}
\renewbibmacro{in:}{}

\newcommand{\Ts}{\mathbf{T}}
\newcommand{\Ss}{\mathbf{S}}
\newcommand{\Ds}{\mathbf{D}}
\newcommand{\Ms}{M}
\newcommand{\Tab}{\Ts_{-\infty,\beta}}
\newcommand{\Taz}{\Ts_{-\infty,0}}
\newcommand{\Tbt}{\Ts_\beta(\theta)}
\newcommand{\Tzt}{\Ts_{0}(\theta)}
\newcommand{\Kop}{\mathbf{K}}
\newcommand{\Pab}{\mathbf{P}_{-\infty}}
\newcommand{\Sb}{\Ss_\beta}
\newcommand{\Sz}{\Ss_{0}}
\newcommand{\Sin}{\Ss_{(n)}}
\newcommand{\Sinm}{\Ss_{(n-1)}}
\newcommand{\Siz}{\Ss_{(0)}}
\newcommand{\Sio}{\Ss_{(1)}}
\newcommand{\Sab}{\Ss_{-\infty}}
\newcommand{\Sabn}{\Ss_{-\infty,(n)}}
\newcommand{\Sop}{\Ss}
\newcommand{\Aop}{\mathbf{A}}
\newcommand{\Bop}{\mathbf{B}}
\newcommand{\Fop}{\mathbf{F}}
\newcommand{\Dn}{\Ds_{n}}
\newcommand{\Dj}{\Ds_{j}}
\newcommand{\PP}{\mathbb{P}}
\newcommand{\Mt}{\Ms_{\tau}}
\newcommand{\Mtt}{\Ms_{\tau(\theta)}}

\title{
  Leray--Hopf Nonuniqueness for Forced 3D Navier--Stokes with Logarithmically
  Supercritical Diffusion at the Lions Exponent
}

\author{Calvin Khor and Xiaoyan Su}

\date{\today}

\begin{document}

\maketitle

\begin{abstract}
  We construct two distinct Leray--Hopf solutions, with the same force and zero
  initial datum, for the three-dimensional forced Navier--Stokes equations with
  generalised diffusion \(m( |\nabla| )\) given by the symbol
  \(m(r) = r^{2\alpha}\ell(r)\), where \(\alpha = \frac{5}{4}\) is the Lions
  exponent, \(\ell\) is positive, smooth, and smoothly varying with index zero,
  and \(\ell \in L^2\left( (1,\infty);\frac{dr}{r} \right)\). The force belongs
  to the energy-dual space \(L_t^2\mathcal{V}'_m\) and is smooth for positive
  times. We also show that this condition is essentially optimal, with
  uniqueness when \(\ell\) is monotone decreasing to 0 and
  \(\ell \notin L^2\left( (1,\infty);\frac{dr}{r} \right)\). In particular, the
  construction shows nonuniqueness of solutions for every \(q > \frac{1}{2}\) to
  \[
    \partial_{t}u + u \cdot \nabla u + \nabla p + \frac{( - \Delta)^{\frac{5}{4}}}{\left( \log\left( e + ( - \Delta)^{\frac{1}{2}} \right) \right)^q}u
    = f,\quad\nabla \cdot u
    = 0,
  \]
  and we have uniqueness of solutions for $q \le \frac12$. For
  \(\frac12 < \alpha < \frac{5}{4}\), and for the above equation at the Lions
  exponent with \(q > 1\), we also show Leray--Hopf nonuniqueness for
  \(f \in L_t^2\mathcal{V}'_m \cap L_t^{1}L_x^2\).
\end{abstract}

\section{Introduction}
\label{sec:introduction}

In this paper we consider the well-posedness of Leray--Hopf solutions to the
forced three-dimensional Navier--Stokes equations with generalised diffusion,
\begin{equation}
  \partial_{t}u + u \cdot \nabla u + \nabla p + m( |\nabla| )u
  = f,\quad\nabla \cdot u
  = 0.\label{eq:main}
\end{equation}
Here \(x \in \mathbb{R}^3\) is position, \(t\) is time, \(u\), \(p\), and \(f\)
are velocity, pressure, and force respectively, and \(m( |\nabla| )\) denotes
the radial Fourier multiplier defined by
\((m( |\nabla| )v)^\wedge(\zeta) = m( |\zeta| )\hat{v}(\zeta).\) The classical
Navier--Stokes equations,
\begin{equation}
  \partial_{t}u + u \cdot \nabla u + \nabla p - \Delta u = f,\quad\nabla \cdot u
                                                         = 0,\label{eq:main-ns}
\end{equation}
correspond to \(m(r) = r^2\), and the fractional forced Navier--Stokes
equations,
\begin{equation}
  \partial_{t}u + u \cdot \nabla u + \nabla p + ( - \Delta)^{\alpha}u
  = f,\quad\nabla \cdot u
  = 0,\label{eq:fractional-ns}
\end{equation}
which generalise~\eqref{eq:main-ns}, correspond to the power symbols
\(m(r) = r^{2\alpha}\), \(\alpha > 0\).

It has in recent years become clear that the Navier--Stokes
equations~\eqref{eq:main-ns} are ill-behaved in various senses, for instance
there is (with $f = 0$) norm inflation in \(\dot B^{-1}_{\infty,\infty}\) by
Bourgain and Pavlovi\'c\cite{bourgain-pavlovic}, the nonuniqueness of
distributional solutions by Buckmaster and Vicol~\cite{buckmaster-vicol}, the
nonuniqueness of solutions for data in $\text{BMO}^{-1}$ by Coiculescu and
Palasek~\cite{coiculescu-palasek}, and for data in
$C^\infty(\mathbb R^3 \setminus \{0\}) \cap L^q$ by Hou, Wang, and
Yang~\cite{hou-wang-yang} with a computer assisted proof. There is also the more
recent work~\cite{ionescu-jia-palasek-1, cadiot} of Ionescu, Jia, and Palasek
which prove nonuniqueness in the unforced, axially symmetric swirl-free case,
also with computer assistance.

For Leray--Hopf solutions to~\eqref{eq:main-ns} with say, \(C_c^\infty\) forcing
and initial data, the nonuniqueness of solutions requires a loss of regularity,
as the solution coincides with the classical smooth solution for short times.
Very recently, it was shown that Leray--Hopf solutions can blow up in finite
time starting from smooth data, at first with a forcing
term~\cite{OpenAI2026NavierStokes}, and then even with $f = 0$
\cite{alpoge-buckmaster-ns}. There was also blowup
established~\cite{OpenAI2026Euler, alpoge-buckmaster-forced-euler,
alpoge-buckmaster-coiculescu-euler} for the related Euler equations, also with
smooth initial data, and without forcing.

The result of this paper, which is a generalisation of the work of Albritton,
Brué, and Colombo~\cite{abc} and then of the current authors' previous
collaboration with Miao~\cite{kms}, proves instantaneous nonuniqueness after the
initial time while remaining inside the finite energy class, using an
instability discovered by Vishik~\cite{vishik-1, vishik-2}. Like the earlier
mentioned works~\cite{hou-wang-yang, ionescu-jia-palasek-1}, these results serve
to illustrate the unintuitive fact that the strict questions of nonuniqueness
and persistence of regularity for Leray--Hopf solutions can be discussed
somewhat independently. Indeed, the recent preprint of Galdi and
Gazzola~\cite{galdi-gazzola-2026} gives solutions that exhibit both blowup and
uniqueness simultaneously. We mention in passing the exposition~\cite{abcdgjk}
and the subsequent simplified proof~\cite{zbMATH08060557} of Vishik's result.

We remark here further advancements in this direction for the interested reader.
Albritton, Bru\'e and Colombo extended~\cite{gluing} the result to the torus and
bounded domains. Hofmanov\'a, Zhu, and Zhu~\cite{hofmanova-zhu-zhu-forced}
showed that forces producing Leray--Hopf nonuniqueness from zero initial data
are dense in \(L_t^1L_x^2\), and that any prescribed divergence-free \(L^2\)
initial datum admits a force producing nonuniqueness. There are also results for
stochastic forced equations. Hofmanov\'a, Zhu, and
Zhu~\cite{hofmanova-zhu-zhu-forced} also extended the instability construction
to linear multiplicative noise, while Bru\'e, Jin, Li, and
Zhang~\cite{brue-jin-li-zhang} treated additive noise and established the
sharpness of the Lions exponent for uniqueness in law in the Leray class. Chen,
Dong, and Zheng~\cite{chen-dong-zheng} obtained nonuniqueness below the Lions
exponent for linear multiplicative noise, including local nonuniqueness with a
deterministic force, robustness under suitable perturbations of this force, and
global nonuniqueness with a suitably chosen random force.

Lions showed~\cite{lions} that there is global well-posedness of energy class
solutions to~\eqref{eq:fractional-ns} for \(f = 0\) and
\(\alpha \geq \frac{5}{4}\). It was later shown by the current authors in joint
work with Miao~\cite{kms} that for \(\frac{1}{2} < \alpha < \frac{5}{4}\) there
is nonuniqueness of energy class solutions to the fractional Navier--Stokes
equations with \(f \in L^1_tL^2_x\), which shows the optimality of the Lions
exponent \(\alpha = \frac{5}{4}\) for the forced fractional Navier--Stokes
equations.

Having settled this question for fractional power multipliers
\(m(\rho) = \rho^{2\alpha}\), we consider lower order improvements of this
result at the critical Lions exponent \(\alpha = \frac{5}{4}\). The
well-posedness above the Lions exponent was refined to the generalised
equation~\eqref{eq:main} first by Tao~\cite{tao-log} and later in the optimal
form by Barbato--Morandin--Romito~\cite{barbato-morandin-romito-apde}. Tao
showed global well-posedness assuming that for $\rho \gg 1$ we have (in our
notation) \(m(\rho) \geq \rho^{\frac{5}{2}}\ell \geq 0\) with $\ell$ eventually
nonincreasing, and satisfies
\begin{equation}
  \ell \notin L^2\left( (1,\infty);\frac{dr}{r} \right)\label{eq:tao-bound}
\end{equation}
He then conjectured that it suffices to replace~\eqref{eq:tao-bound} with
\begin{equation}
  \ell \notin L^1\left( (1,\infty);\frac{dr}{r} \right)\label{eq:romito-bound}
\end{equation}
which Barbato--Morandin--Romito proved (under the additional assumption that
$\rho^{5/2}\ell$ is eventually nondecreasing). For the unforced problem, their
result gives global regularity for the logarithmically supercritical multiplier
\(m_{q}(r) = r^{\frac{5}{2}}( \log(e + r) )^{- q}\) when \(q \leq 1\). Recently,
Romito and Triggiano showed~\cite{romito-triggiano} that, for \(q > 1\), there
is nonuniqueness of distributional solutions to the unforced equation with
\(m_{q}( |\nabla| )\), which proves in a sense the sharpness of the exponent,
and gives evidence for the sharpness of~\eqref{eq:romito-bound} in general.

In this paper, we show that for Leray--Hopf solutions to the generalised
Navier--Stokes equations~\eqref{eq:main} with forcing in the natural energy-dual
space \(L^2\mathcal{V}'_m\) (\(\mathcal{V}_m\) defined below), somewhat
curiously, it is in fact Tao's condition~\eqref{eq:tao-bound} that is sharp. We
also prove nonuniqueness when \(\ell\) does not satisfy~\eqref{eq:romito-bound},
i.e. $\ell \in L^1(1,\infty; \frac{dr}r)$ for \(f \in L^1_tL^2_x\). Notably,
this leaves open the question of uniqueness for forces \(f \in L^1_tL^2_x\) in
the intermediate regime
\[
  \ell \in L^2 \setminus L^1\left( (1,\infty);\frac{dr}r \right).
\]
For $\ell(r) = (\log(e + r))^{-q}$, this is the range $1/2 < q \le 1$, as
illustrated in Figure~\ref{fig:log-thresholds}.

\subsection{Notation and Main results}

We use the unitary Fourier transform
$\hat v(\zeta) = (2\pi)^{-3/2} \int_{\mathbb R^3}e^{-ix \cdot \zeta}v(x)\,dx$
(the symbol $\xi$ is reserved for our new coordinates below.) For
\(s \in \mathbb{R}\), let \(H_\sigma^{s}( \mathbb{R}^3 )\) be the closed
subspace of \(H^{s}( \mathbb{R}^{3};\mathbb{R}^3 )\) consisting of
divergence-free fields, and put \(L_\sigma^2 \coloneqq H_\sigma^0\). The Leray
projector \(\PP\) is the Fourier multiplier given by the symbol
\(I - \frac{\zeta\zeta^T}{|\zeta|^2}\).

We recall the following definition from the theory of regularly varying
functions~\cite{bingham-goldie-teugels}.

\begin{definition}

  For \(\beta \in {\mathbb{R}}\), a non-negative radial symbol
  \(m:\left[ 0,\infty \right) \to \left[ 0,\infty \right)\) belongs to
  \(\mathrm{SR}_\beta\) if for all \(r \gg 1\),
  \[
    m(r) = r^\beta\ell(r),
  \]
  where \(\ell\) is positive, eventually smooth, and smoothly varying with index
  zero, i.e.
  \begin{gather}
    ( r\partial_r )^k\log\ell(r) \xrightarrow[r \to \infty]{} 0 \quad (k = 1,2,\dots). \label{eq:smooth-slow-weight}
  \end{gather}
\end{definition}
We will additionally impose for our convenience that \(\ell\) is actually smooth
up to \(r = 0\). Although our focus is primarily at the Lions exponent
$\beta = 2\alpha = 5/2$ we will proceed with a general index as the proof allows
for clarity. We use symbols in $\mathrm{SR}_{2\alpha}$ so that the rescaled
diffusion below converges to \(\rho^{2\alpha}\) locally.

For any symbol \(m\), define
\begin{align*}
  {\dot{H}}_{m} & \coloneqq \left\{ v \in L^{2}( \mathbb{R}^3 ):\int_{\mathbb{R}^3}m( |\zeta| )| {\hat{v}(\zeta)} |^{2}d\zeta < \infty \right\},\quad & \| v \|_{{\dot{H}}_m}^2   & \coloneqq \int_{\mathbb{R}^3}m( |\zeta| )| {\hat{v}(\zeta)} |^{2}d\zeta, \\
  \mathcal{V}_m & \coloneqq L_\sigma^{2}( \mathbb{R}^3 ) \cap {\dot{H}}_m ,\quad                                                                      & \| v \|_{\mathcal{V}_m}^2 & \coloneqq \| v \|_{L^2}^2 + \| v \|_{{\dot{H}}_m}^{2}.
\end{align*}
For \(m(r) = r^{2\alpha}\) this is the usual homogeneous \(H^\alpha\) seminorm
on \(L^2\). For an interval \(I\) and a Banach space \(X\), \(C_{w}(I;X)\)
denotes the space of weakly continuous maps from \(I\) to \(X\).

\begin{definition}

  Let \(m:\left[ 0,\infty \right) \to [0,\infty)\),
  \(f \in L^1(0,T;L^2) + L^{2}( 0,T;\mathcal{V}'_m )\), and
  \(u_0 \in L_\sigma^{2}( \mathbb{R}^3 )\). A Leray--Hopf solution
  of~\eqref{eq:main} on \(\lbrack 0,T\rbrack\) is a divergence-free field
  \[
    u \in C_{w}( \lbrack 0,T\rbrack;L_\sigma^2 )
          \cap L^\infty( 0,T;L_\sigma^2 )
          \cap L^{2}( 0,T;\mathcal{V}_m )
  \]
  which satisfies the equation distributionally, attains \(u_0\) strongly in
  \(L^2\) as \(t \downarrow 0\), and for which there is a full-measure set
  \(\mathcal{G} \subset \lbrack 0,T\rbrack\), containing \(0\), such that
  \begin{equation}
    \frac{1}{2}\| {u(t)} \|_{L^2}^2 + \int_s^t\| {u(t')} \|_{{\dot{H}}_m}^2\, dt'
    \leq \frac{1}{2}\| {u(s)} \|_{L^2}^2
         + \int_s^t\langle f(t'),u(t')\rangle\, dt'\label{eq:forced-energy}
  \end{equation}
  for every \(s \in \mathcal{G}\) and every \(t \in \lbrack s,T\rbrack\). The
  last term denotes the dual pairing.
\end{definition}

\begin{theorem}
  \label{thm:main}

  Let \(\frac{1}{2} < \alpha \leq \frac{5}{4}\) and let
  \(m = r^{2\alpha}\ell \in \mathrm{SR}_{2\alpha}\) with
  \(\ell \in C^\infty\left( \left[ 0,\infty \right),(0,\infty) \right)\).
  Suppose also in the case \(\alpha = \frac{5}{4}\) that
  \(\ell \in L^2\left( (1,\infty);\frac{dr}{r} \right)\). Then there exist a
  time \(t_\ast > 0\), a real-valued divergence-free force
  \[
    f \in L^{2}( 0,t_\ast;\mathcal{V}'_m )
          \cap C_{\text{loc}}^\infty( (0,t_\ast] \times \mathbb{R}^3 ),
  \]
  and two distinct real-valued Leray--Hopf solutions \(u^{(1)},u^{(2)}\)
  to~\eqref{eq:main} on \(\lbrack 0,t_\ast \rbrack\), both
  satisfying~\eqref{eq:forced-energy} with equality and
  \(\mathcal{G} = \lbrack 0,t_\ast \rbrack\), and
  \[
    u^{(j)}(t)
    \to        0\quad\text{ strongly in }L^{2}( \mathbb{R}^3 )\quad\text{ as }t
    \downarrow 0.
  \]
  Moreover, if \(\alpha < \frac{5}{4}\), or if \(\alpha = \frac{5}{4}\) and
  \(\ell \in L^1\left( (1,\infty);\frac{dr}{r} \right)\) the same construction
  gives
  \[
    f \in L^{2}( 0,t_\ast;\mathcal{V}'_m )
          \cap C_{\text{loc}}^\infty( (0,t_\ast] \times \mathbb{R}^3 )
          \cap L^{1}( 0,t_\ast;L^{2}( \mathbb{R}^3 ) ).
  \]
\end{theorem}


To show that this result is sharp, we have the following.

\begin{theorem}
  \label{thm:local-wellposedness}
  Let \(\alpha \geq \frac{5}{4}\) and let
  \(m = r^{2\alpha}\ell \in \mathrm{SR}_{2\alpha}\) with
  \(\ell \in C^\infty\left( \left[ 0,\infty \right),(0,\infty) \right)\) and let
  \(\underline{\ell}(r) \coloneqq \min(1,\min\limits_{1 \leq r' \leq r}\ell(r'))\)
  for \(r \geq 1\). If \(\alpha > \frac{5}{4}\), or if \(\alpha = \frac{5}{4}\)
  and
  \[
    \int_1^\infty{\underline{\ell}(r)}^2\frac{dr}{r} = \infty,
  \]
  then any two Leray--Hopf solutions \(u^{(1)},u^{(2)}\) of~\eqref{eq:main} with
  the same \(L^2\) initial data and the same forcing
  \(f \in L^{1}L^2 + L^2\mathcal{V}'_m\) must coincide.
\end{theorem}
From the proof, it suffices that the \(u^{(i)}\) are distributional solutions in
$C_{w}( \lbrack 0,T\rbrack;L_\sigma^2 ) \cap L^\infty( 0,T;L_\sigma^2 ) \cap L^{2}( 0,T;\mathcal{V}_m ).$
The need to use the capped lower envelope \(\underline{\ell}\) matches with the
lower bounding of the symbol in~\cite{tao-log}
and~\cite{barbato-morandin-romito-apde}. Of course, if \(\ell\) is monotone
decreasing to 0 then \(\underline{\ell} = \ell\) eventually.

For the logarithmically modified multiplier
\(m_q(r) = \frac{r^{5/2}}{(\log(e + r))^q}\), we have the following corollary.
\begin{corollary}
  \label{cor:log}
  If \(q \le \frac12\), then any two distributional solutions of
  \begin{align}
    \label{eq:log}
    \partial_{t}u + u \cdot \nabla u + \nabla p + \frac{( - \Delta)^{\frac{5}{4}}}{\left( \log\left( e + ( - \Delta)^{\frac{1}{2}} \right) \right)^q}u
    = f,\quad\nabla \cdot u
    = 0.
  \end{align}
  in
  \(C_w([0,T];L_\sigma^2) \cap L^\infty(0,T;L_\sigma^2) \cap L^2(0,T;\mathcal V_{m_q})\)
  with the same initial data and the same force
  \(f \in L^1(0,T;L^2) + L^2(0,T;\mathcal V_{m_q}')\) coincide.

  If $q > \frac12$, there exist $t_\ast > 0$, a force
  \(f \in L^2(0,t_\ast;\mathcal V_{m_q}') \cap C^\infty_{\mathrm{loc}}((0,t_\ast] \times \mathbb R^3),\)
  and two distinct Leray--Hopf solutions to~\eqref{eq:log} with zero initial
  datum which satisfy the energy equality. If $q > 1$, the force can be chosen
  in
  $L^2(0,t_\ast;\mathcal V_{m_q}') \cap C^\infty_{\mathrm{loc}}((0,t_\ast] \times \mathbb R^3) \cap L^1(0,t_\ast;L^2)$.
\end{corollary}

\definecolor{myred}{HTML}{ee4136}
\definecolor{myblue}{HTML}{0074d9}
\definecolor{myyellow}{HTML}{eedc00}
\begin{figure}[htbp]
  \centering
  \begin{tikzpicture}[
      every node/.style={font=\small, align=center},
      uniqueness/.style={ draw=myblue, very thick, -{Stealth} },
      nonuniqueness/.style={ draw=myred, very thick, -{Stealth} },
      openregion/.style={ draw=myyellow, very thick, dashed },
      guide/.style={gray!60, densely dotted}
    ]

    \def\qleft{0}
    \def\qhalf{3.5}
    \def\qone{7}
    \def\qright{13}

    \def\yunforced{2.7}
    \def\ydual{1.8}
    \def\yone{0.5}


    \draw[guide] (\qhalf,0) -- (\qhalf,3);
    \draw[guide] (\qone,0) -- (\qone,3);

    \draw[Stealth-Stealth] (\qleft,0) -- (\qright,0) node[right] {$q$};

    \draw (\qhalf,-0.08) -- (\qhalf,0.08) node[below=5pt] {$\frac12$};
    \draw (\qone,-0.08) -- (\qone,0.08) node[below=5pt] {$1$};


    \node[anchor=east, align=right] at (-0.25,\yunforced) {$f = 0$};

    \draw[uniqueness] (\qone,\yunforced) -- node[above=3pt]
      { {\color{myblue}global regularity}~\cite{barbato-morandin-romito-apde} }
      (\qleft,\yunforced);

    \draw[nonuniqueness] (\qone,\yunforced) -- node[above=3pt, xshift=0pt]
      { {\color{myred}distributional nonuniqueness}~\cite{romito-triggiano} }
      (\qright,\yunforced);


    \node[anchor=east, align=right] at (-0.25,\ydual)
      {$f \in L_t^2\mathcal{V}_{m_q}'$};

    \draw[uniqueness] (\qhalf,\ydual) -- node[above=3pt] {} (\qleft,\ydual);

    \draw[nonuniqueness] (\qhalf,\ydual) -- node[above=3pt,xshift=0pt] {  }
      (\qright,\ydual);


    \node[anchor=east, align=right] at (-0.25,\yone) {$f \in L_t^1L_x^2$};

    \draw[uniqueness] (\qhalf,\yone) -- node[above=3pt]
      {\color{myblue}uniqueness \\
        \color{myblue}(Theorem~\ref{thm:local-wellposedness})}
      (\qleft,\yone);

    \draw[nonuniqueness] (\qone,\yone) -- node[above=3pt]
      {
        \color{myred} Leray--Hopf nonuniqueness \\
        \color{myred}(Theorem~\ref{thm:main})
      }
      (\qright,\yone);


    \draw[myblue, thick, fill] (\qhalf,\yone) circle (2pt);
    \draw[myblue, thick, fill] (\qhalf,\ydual) circle (2pt);
    \draw[myblue, thick, fill] (\qone,\yunforced) circle (2pt);

    \draw[myred,thick, fill=white] (\qone,\yone) circle (2pt);

  \end{tikzpicture}
  \caption{
    An illustration of uniqueness and nonuniqueness results for the
    logarithmically modified multiplier $m_q(r) = r^{5/2}(\log(e + r))^{-q}$.
    The unforced results are due to
    Barbato--Morandin--Romito~\cite{barbato-morandin-romito-apde} and
    Romito--Triggiano~\cite{romito-triggiano}. The results with $f \neq 0$ are
    stated in Corollary~\ref{cor:log} and follow from Theorems~\ref{thm:main}
    and~\ref{thm:local-wellposedness}. In the range $q \in (1/2,1]$, the
    question of uniqueness of Leray--Hopf solutions with $f \in L_t^1L_x^2$ is
    not covered by the results of this paper.
  }
  \label{fig:log-thresholds}
\end{figure}
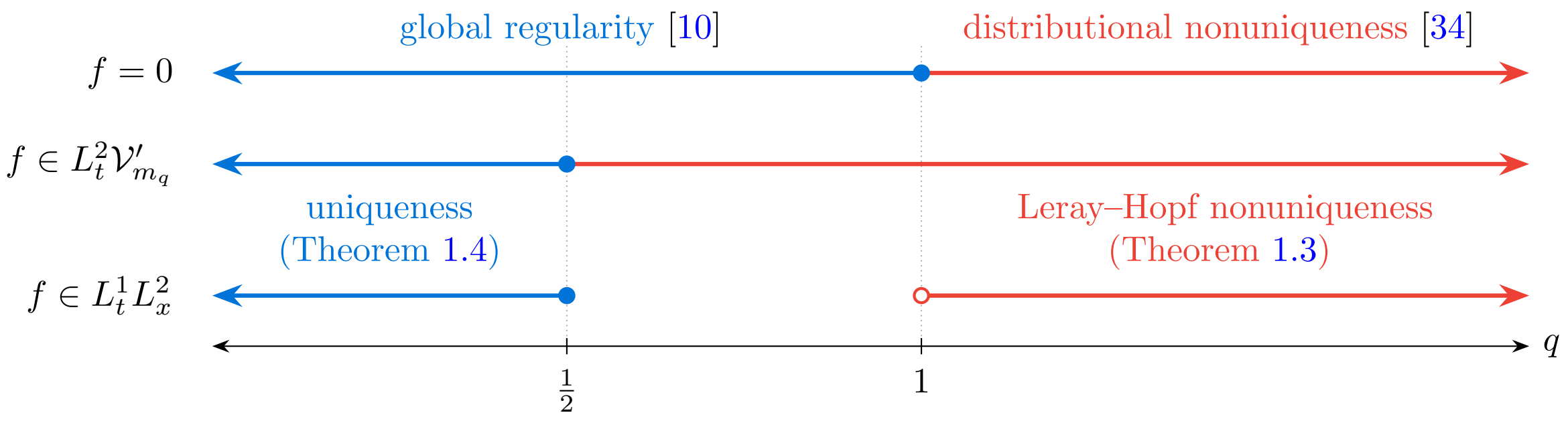

\subsection{Proof method and organisation}

For the nonuniqueness Theorem~\ref{thm:main}, we follow the same strategy of
studying the self-similar equations as in~\cite{abc} and~\cite{kms}, but the
details change because the multiplier is not homogeneous. If we were to use the
same scaling as in~\cite{kms}, then at the Lions endpoint, the diffusion term
would converge to zero, a problem if we want to use parabolic theory, and our
solutions would not converge to 0 in $L^2$ as \(t \downarrow 0\). Our choice of
modified coordinates introduced in Section~\ref{sec:2} below avoids both issues.
This means that we only have almost self-similar variables. The good news is
that the rescaled diffusion \(\Mt( |\nabla| )\) behaves asymptotically like the
fractional Laplacian \(( - \Delta)^\alpha\) as we approach logarithmic time
\(\tau = - \infty\), so the linear instability of~\cite{kms} is still present.
In Section~\ref{sec:ancient-instability}, we recall the various results
of~\cite{kms} that we need.

The analysis is made slightly more technical by the fact that the rescaled
diffusion \(\Mt( |\nabla| )\) is not independent of (logarithmic) time \(\tau\),
which means we need to use some constructions for the theory of non-autonomous
operators. This is encapsulated in Proposition~\ref{prop:ancient-unstable-mode},
which Section~\ref{sec:ancient-mode} is devoted to proving. Indeed, it should be
noted that we could have used more abstract theory of nonautonomous evolution
(in particular the interested reader can consult~\cite{rrsv-2000,
backes-dragicevic} and the references therein), but the presentation below is
self-contained and avoids tedious checking of hypotheses such as sectoriality.
Section~\ref{sec:construction} constructs the two solutions with Golovkin's
trick~\cite{golovkin-1964, golovkin-1965} and checks the assertions of
Theorem~\ref{thm:main}: convergence to 0 initial data, dissipation, force
estimates, and the energy identity.

We remark that we could have instead pursued a nonuniqueness proof with a
quasilinear flavour using the ideas of~\cite{ac} which chooses a similarity
scale on which the diffusion becomes perturbative in the setting of the 2D
hypodissipative Navier--Stokes equations. One would expect such a proof to allow
also $\alpha < 1/2$, but as we are interested in the wellposedness near the
Lions exponent, we preferred to keep the diffusion operator in the forefront in
our exposition to directly exploit the linear instability of~\cite{kms} as
described above.

In Section~\ref{sec:uniqueness}, we prove the uniqueness result by estimating
the difference of two solutions in \(H^{- 1}\) and applying the Osgood
uniqueness lemma. The nonlinear estimate in the space $\mathcal V_m$ which is
adapted to the multiplier \(m\) (Lemma~\ref{lem:uniqueness-nonlinear-est})
allows us to close the argument.

\section{The equations in almost self-similar variables}
\label{sec:2}

The following lemma concerns basic properties of our multipliers and motivates
the change of variables.

\begin{lemma}
  \label{lem:symbol-calculus}

  Let \(m = r^{2\alpha}\ell \in \mathrm{SR}_{2\alpha}\) with
  \(\ell \in C^\infty\left( \left[ 0,\infty \right),(0,\infty) \right)\) and
  \(\alpha > \frac{1}{2}\). Then \(m\) is strictly increasing and unbounded for
  all sufficiently large arguments. In addition, for every integer \(N \geq 0\),
  \begin{equation}
    \frac{m(R\rho)}{m(R)} \to \rho^{2\alpha}\quad\text{ in }\quad C_{\text{loc}}^N\left( (0,\infty)_\rho \right)\quad\text{ as }\quad R
                          \to \infty.\label{eq:regular-variation}
  \end{equation}
  Moreover, for every \(\kappa > 0\) and integer \(j \geq 0\), there are
  constants \(R_\kappa,C_\kappa,C_{j,\kappa} > 0\) such that, whenever
  \(R \geq R_\kappa\) and \(\rho > 0\),
  \begin{align}
    \frac{1}{C_\kappa}\rho^{2\alpha}\min(\rho^\kappa,\rho^{- \kappa}) & \leq \frac{m(R\rho)}{m(R)} \leq C_\kappa\rho^{2\alpha}\max(\rho^\kappa,\rho^{- \kappa}),\label{eq:potter-bounds}
  \end{align}

  \begin{align}
    \left| {( \rho\partial_\rho )^j\frac{m(R\rho)}{m(R)}} \right| & \leq C_{j,\kappa}\rho^{2\alpha}\max(\rho^\kappa,\rho^{- \kappa}).\label{eq:potter-bounds-2}
  \end{align}
  Finally, if either $\frac12 < \alpha < \frac54$, or $\alpha = \frac54$ and
  $\ell \in L^2((1,\infty);dr/r)$, then
  \begin{equation}
    r^{-5/2}m(r) \longrightarrow 0 \qquad\text{as }r
                 \to             \infty. \label{eq:symbol-little-o}
  \end{equation}
  If $\frac12 < \alpha < \frac54$, then moreover
  \begin{equation}
    r^{-5/2}m(r) \in L^p((1,\infty);dr/r) \qquad\text{for every }1
                 \le p
                 <   \infty. \label{eq:symbol-Lp}
  \end{equation}
\end{lemma}
\begin{proof}
  Since
  \begin{equation}
    r\partial_r\log m(r)
    = 2\alpha
      + r\partial_r\log\ell(r) \xrightarrow[r\to\infty]{} 2\alpha,\label{eq:log-derivative}
  \end{equation}
  there is \(R_0 > 0\) so that \(r\partial_r\log m(r) \geq \alpha\) for
  \(r \geq R_0\). Thus \(m\) is strictly increasing on this tail, and
  \[
    \log\frac{m(r)}{m( R_0 )}
    =    \int_{R_0}^{r}s\partial_s\log m(s)\frac{ds}{s}
    \geq \alpha\log\frac{r}{R_0} \xrightarrow[r \to \infty]{} \infty.
  \]
  Hence \(m\) is unbounded.

  Next we prove~\eqref{eq:regular-variation}. For fixed \(\rho > 0\),
  \begin{align*}
    \log\frac{m(R\rho)}{m(R)\rho^{2\alpha}} = \log\ell(R\rho) - \log\ell(R) & = \int_1^{\rho}s\partial_{s}( \log\ell(Rs) )\frac{ds}{s}                                                        \\
                                                                            & = \int_1^\rho\sigma\partial_\sigma( \log\ell(\sigma) )|_{\sigma=Rs}\frac{ds}{s} \xrightarrow[R \to \infty]{} 0.
  \end{align*}
  This gives the locally uniform convergence
  \(\frac{m(R\rho)}{m(R)} \to \rho^{2\alpha}\) in \(C^0\). For the general case,
  recall the identity
  \(( \rho\partial_\rho )^{j}( f(R\rho) ) = ( x\partial_x )^{j}f(x) |_{x=R\rho}\).
  Thus, we have for \(k \geq 1\) that
  \[
    ( \rho\partial_\rho )^k\log\frac{m(R\rho)}{m(R)\rho^{2\alpha}}
    = ( \rho\partial_\rho )^k\log\ell(R\rho)
    = ( r\partial_r )^k\log\ell(r) |_{r=R\rho} \xrightarrow[R \to \infty]{} 0
  \]
  from which the locally uniform convergence in \(C^N\) follows.

  The scalar estimate in~\eqref{eq:potter-bounds} follows
  from~\cite[Theorem 1.5.6]{bingham-goldie-teugels}. Positivity and smoothness
  of \(\ell\) on \(\left[ 0,\infty \right)\) extend it uniformly to the
  low-frequency range \(R\rho < R_\kappa\).

  For \(j \geq 1\), the chain rule gives
  \[
    ( \rho\partial_\rho )^j\frac{m(R\rho)}{m(R)} = \frac{m(R\rho)}{m(R)}\left. \frac{( r\partial_r )^{j}m(r)}{m(r)} \right|_{r=R\rho}.
  \]
  It is easy to see that \(\frac{( r\partial_r )^{j}m(r)}{m(r)}\) is a
  polynomial in \(( r\partial_r )\log m(r),\ldots,( r\partial_r )^j\log m(r)\)
  with coefficients depending only on \(j\). These logarithmic derivatives are
  bounded near zero because \(\ell\) is positive and smooth there, and they are
  bounded at infinity by~\eqref{eq:smooth-slow-weight}. Thus,
  \(\frac{( r\partial_r )^{j}m(r)}{m(r)}\) is bounded on \((0,\infty)\) for
  every \(j \geq 1\). Together with~\eqref{eq:potter-bounds}, this
  proves~\eqref{eq:potter-bounds-2}.

  If \(\alpha < \frac{5}{4}\), choose \(\kappa \ll 1\) and put
  \(\varepsilon = \frac{5}{2} - 2\alpha - \kappa > 0\). The upper bound
  in~\eqref{eq:potter-bounds} gives
  \[
    r^{- \frac{5}{2}}m(r) \leq Cr^{- \varepsilon},\quad r \geq R_\kappa.
  \]
  This proves both~\eqref{eq:symbol-little-o} and~\eqref{eq:symbol-Lp}.

  It remains to prove~\eqref{eq:symbol-little-o} when \(\alpha = \frac{5}{4}\).
  Put \(g(x) = \ell(e^{x})\). Then \(g \in L^2\left( (0,\infty) \right)\) and
  \(( \log g )'(x) \to 0\). If \(g\) did not tend to zero, there would exist
  \(x_n \to \infty\) with \(g( x_n ) \geq \varepsilon > 0\). The intervals
  \(\lbrack x_n - 1,x_n + 1 \rbrack\) may be taken disjoint and
  \(g(x) \geq e^{- 1}\varepsilon\) on each of them, contradicting
  \(g \in L^2\left( (0,\infty) \right)\). Hence \(\ell(r) \to 0\), as needed.
\end{proof}

\begin{remark}
  In particular, every symbol satisfying the hypotheses of
  Theorem~\ref{thm:main} satisfies
  \begin{equation*}
    r^{- \frac{5}{2}}m(r) \in L^2\left( (1,\infty),\frac{dr}{r} \right).
  \end{equation*}
  For \(\alpha < \frac{5}{4}\), this is part of Lemma~\ref{lem:symbol-calculus}
  at \(p = 2\); at the Lions endpoint, \(r^{- \frac{5}{2}}m(r) = \ell(r)\), so
  we need to impose it as an assumption in Theorem~\ref{thm:main}. Later in
  Section~\ref{sec:construction}, this condition gives the finite dissipation
  in~\eqref{eq:construction-dissipation}, and the energy-dual force estimate.
\end{remark}

\begin{remark}
  \label{rem:multiplier-convergence}
  The above imply that for \(v \in H^{2\alpha+\varepsilon}\),
  \(\frac{m( R|\nabla| )}{m(R)}v \to ( - \Delta)^{\alpha}v\) in $L^2$ as
  \(R \to \infty\).
\end{remark}
By~\eqref{eq:log-derivative}, \(m \to \infty\) as $r \to \infty$ and is
eventually strictly increasing. For all sufficiently small \(t > 0\), let
\(R(t)\) be the unique solution of \(m( R(t) ) = t^{- 1}\). Then
\(R(t) \to \infty\) as \(t \downarrow 0\), and \(R\) is \(C^\infty\) by the
inverse function theorem.

Define the logarithmic time \(\tau = \log t\), frequency scale
\(R_\tau \coloneqq R( e^\tau ) = R(t)\), and
\begin{equation}
  \Mt(\rho)
  \coloneqq \frac{m( R_\tau\rho )}{m( R_\tau )}
  =         e^{\tau}m( R_\tau\rho ),\quad\rho
  \geq      0,\mspace{27mu} b(\tau)
  \coloneqq - \frac{d\log R_\tau}{d\tau}
  =         \frac{m( R_\tau )}{R_{\tau}m'( R_\tau )}.\label{eq:Mtau-b-def}
\end{equation}
The displayed formula for \(b\) follows from differentiating
\(\log m( R_\tau ) = - \tau\). Since
\(R_\tau \to \infty\),~\eqref{eq:log-derivative} gives
\[
  b(\tau) \to \frac{1}{2\alpha}\quad\text{ as }\quad\tau \to - \infty.
\]
Equation~\eqref{eq:regular-variation} identifies the limiting symbol
in~\eqref{eq:Mtau-b-def} as \(\rho^{2\alpha}\). Set
\begin{align*}
  \xi    & = R(t)x,                            & u(x,t)       & = \frac{1}{tR(t)}U(\xi,\tau),    \\
  p(x,t) & = \frac{1}{t^2{R(t)}^2}P(\xi,\tau), & \quad f(x,t) & = \frac{1}{t^{2}R(t)}F(\xi,\tau) 
\end{align*}
Here \(\xi\) is the rescaled space variable, \(U\) is the rescaled velocity,
\(P\) is the rescaled pressure, and \(F\) is the rescaled force. Direct
substitution into~\eqref{eq:main} gives
\begin{equation}
  \begin{split}
    \partial_{\tau}U + ( b(\tau) - 1 - b(\tau)\xi \cdot \nabla_\xi )U + U \cdot \nabla_{\xi}U + \nabla_{\xi}P + \Mt( |\nabla| )U & = F,\mspace{27mu} \\
    \nabla_\xi \cdot U                                                                                                           & = 0.
  \end{split}
  \label{eq:similarity}
\end{equation}

\section{Linear instability at large negative logarithmic time}
\label{sec:ancient-instability}

The appearance of \(\Mt( |\nabla| )\) in~\eqref{eq:similarity} which converges
to the standard fractional Laplacian \(( - \Delta)^\alpha\) as
\(\tau \to - \infty\) suggests to us that we can transfer the instability for
that equation which was already studied in~\cite{kms}. In this section we cite
the pertinent results of that paper.

Let \(\overline{\Omega} \in C_c^\infty( \mathbb{R}^{3};\mathbb{R}^3 )\) be the
axisymmetric pure-swirl unstable Euler vortex used in~\cite{abc} and~\cite{kms},
and let
\(\overline{U} \coloneqq ( - \Delta)^{- 1}\nabla \times \overline{\Omega}\) be
the corresponding velocity. Recall that \(\int\overline{\Omega} = 0\), while
\(\hat{\overline{U}}\) and its radial derivative are bounded at the origin and
decay faster than any power as \(|\zeta| \to \infty\). For
\(\alpha > 0,\ \beta \ge 0\), define the operators
\begin{align*}
  \Kop v & \coloneqq  {\PP} (\overline{U} \cdot \nabla_{\xi}v + v \cdot \nabla_\xi\overline{U}),                                                    \\
  \Tab   & \coloneqq \left( 1 - \frac{1}{2\alpha} \right)I + \left( \frac{1}{2\alpha} \right)\xi \cdot \nabla_\xi - |\nabla|^{2\alpha} - \beta\Kop. 
\end{align*}
As a densely defined operator on \(L_\sigma^2\), we have
\begin{equation*}
  D( \Tab ) = D( \Taz ) = \{ v \in H_\sigma^{2\alpha}( \mathbb{R}^3 ):\xi \cdot \nabla_{\xi}v \in L^{2}( \mathbb{R}^3 ) \},
\end{equation*}
with graph norm
\(\| v \|_{D(\Taz)} \coloneqq \| v \|_{L^2} + \| {\Taz v} \|_{L^2} \cong \| v \|_{L^2} + \| {|\nabla|^{2\alpha}v} \|_{L^2} + \| {\xi \cdot \nabla_{\xi}v} \|_{L^2}. 
\)
\begin{lemma}
  \label{lem:endpoint-instability}

  Let \(\alpha \in \left( \frac{1}{2},\frac{5}{4} \right]\). For every
  \(A > 0\), there is \(\beta_0 = \beta_{0}(A,\alpha)\) such that, for
  \(\beta > \beta_0\), the spectrum of \(\Tab\) contains an isolated spectral
  value in \(\{ z:\Re z > A \}\).
\end{lemma}

\begin{proof}

  This is the velocity form of~\cite[Lemmas 2.3--2.5]{kms}. The endpoint is
  included in the linear argument: as noted in~\cite[Remark 1]{kms}, the
  restriction \(\alpha < \frac{5}{4}\) enters only in the final estimate of the
  paper.
\end{proof}

\begin{definition}[Essential semigroup growth]
  \label{def:essential-growth}

  For a bounded operator \(\mathbf{S}\) on a Banach space \(X\), set
  \[
    \| \mathbf S \|_{\text{ess}} \coloneqq \inf\{ \| {\mathbf S - \mathbf K} \|:\mathbf K\text{  is compact on  }X \}.
  \]
  If \(A\) generates the strongly continuous semigroup \(e^{h \mathbf A}\),
  define
  \begin{equation*}
    \omega_{\text{ess }}( \mathbf A) \coloneqq \inf\{ \omega \in \mathbb{R}:\text{ there is  }C_\omega\text{  such that  }\| e^{h \mathbf A} \|_{\text{ess }} \leq C_{\omega}e^{\omega h}\text{  for every  }h \geq 0 \}.
  \end{equation*}
\end{definition}

\begin{lemma}
  \label{lem:ess}
  Let \(\frac{1}{2} < \alpha \leq \frac{5}{4}\). For \(\beta > 0\), \(\Tab\)
  generates a strongly continuous semigroup and
  \begin{equation*}
    \omega_{\text{ess }}( \Tab ) \leq 1 - \frac{5}{4\alpha} \leq 0.
  \end{equation*}
\end{lemma}

\begin{proof}

  This is~\cite[Lemmas 3.3--3.4]{kms}.
\end{proof}

\begin{proposition}
  \label{prop:frozen-dichotomy}

  For \(\beta \gg 1\), there are numbers \(\sigma_{s}, \sigma_c, \sigma_u\)
  satisfying
  \[
    \frac{1}{2} < \sigma_s < \sigma_c < \sigma_u
  \]
  and a nonzero finite-rank spectral projection \(\Pab\) such that, with
  \[
    E_{-\infty} \coloneqq \operatorname{Ran}\Pab,\quad F_{-\infty}
                \coloneqq \operatorname{Ker}\Pab,
  \]
  \begin{enumerate}
    \item both \(E_{-\infty}\) and \(F_{-\infty}\) are invariant, and
          \begin{equation}
            \left\| {e^{h\Tab}( I - \Pab )} \right\|_{L^2 \to L^2}
            \leq Ce^{\sigma_{s}h},\quad h
            \geq 0.\label{frozen-stable}
          \end{equation}
    \item The restriction of \(e^{h\Tab}\) to \(E_{-\infty}\) is invertible,
          with
          \begin{equation}
            \left\| \left( e^{h\Tab}|_{E_{-\infty}} \right)^{- 1} \right\|_{E_{-\infty} \to E_{-\infty}}
            \leq Ce^{- \sigma_{u}h},\quad h
            \geq 0.\label{frozen-unstable}
          \end{equation}
    \item The projection \(\Pab\) preserves real-valued fields.
  \end{enumerate}
\end{proposition}

\begin{proof}

  Choose \(\beta\) so that Lemma~\ref{lem:endpoint-instability} gives spectrum
  in \(\{ z:\operatorname{Re}z > 2 \}\). Since
  \(\omega_{\text{ess }}( \Tab ) \leq 0\) by Lemma~\ref{lem:ess}, the spectral
  decomposition theorem for strongly continuous
  semigroups~\cite[Corollary IV.2.11 and Theorem V.3.1]{engel-nagel} implies
  that the spectrum in $\{z:\operatorname{Re}z > 1/2\}$ consists of finitely
  many isolated eigenvalues of finite algebraic multiplicity.

  Choose \(\sigma_c \in \left( \frac{1}{2},2 \right)\) so that there is no
  spectral value with real part \(\sigma_c\), and let \(\Pab\) be the Riesz
  projection onto the spectral values with real part greater than \(\sigma_c\).
  We may then choose \(\sigma_s\) and \(\sigma_u\) with
  \[
    \frac{1}{2} < \sigma_s < \sigma_c < \sigma_u
  \]
  such that every spectral value on
  \(F_{-\infty} \coloneqq \operatorname{Ker}\Pab\) has real part less than
  \(\sigma_s\), while every spectral value on
  \(E_{-\infty} \coloneqq \operatorname{Ran}\Pab\) has real part greater than
  \(\sigma_u\). The invariance follows because the Riesz projection commutes
  with the semigroup. The semigroup estimate~\eqref{frozen-stable} on the
  complementary space follows from~\cite[Theorem V.3.1]{engel-nagel}, while the
  inverse estimate~\eqref{frozen-unstable} on $E_{-\infty}$ follows from the
  fact that $\dim E_{-\infty} < \infty$. Finally, $\Pab$ preserves real fields
  because $\Tab$ is real, and we can choose the contour defining it to be
  symmetric across the real axis.
\end{proof}

\section{Nonautonomous linear evolution}
\label{sec:ancient-mode}

In this section, we prove the following proposition.

\begin{proposition}
  \label{prop:ancient-unstable-mode}

  Choose \(\beta\) as in Proposition~\ref{prop:frozen-dichotomy}. For all
  \(\tau_\ast \ll - 1\), there exist \(\gamma > \frac{1}{2}\) and a
  divergence-free field \(U^L\) satisfying
  \begin{equation}
    \partial_{\tau}U^L
    = ( 1 - b(\tau) )U^L
      + b(\tau)\xi \cdot \nabla_{\xi}U^L
      - \Mt( |\nabla| )U^L
      - \beta{\PP}(\overline{U} \cdot \nabla_{\xi}U^L + U^L \cdot \nabla_\xi\overline{U})\label{eq:linear-ancient}
  \end{equation}
  distributionally on \(\left( - \infty,\tau_\ast \right] \times \mathbb{R}^3\).
  For every integer \(k \geq 0\),
  \(U^L \in C\left( \left( - \infty,\tau_\ast \right];H_\sigma^k \right)\) and
  \begin{equation}
    \| {U^{L}( \tau_\ast )} \|_{L^2}
    =    1,\quad\| {U^{L}(\tau)} \|_{H^k}
    \leq C_{k}e^{\gamma(\tau-\tau_\ast)}.\label{eq:UL-bounds}
  \end{equation}
  The field \(U^{L}(\tau)\) is nonzero for every \(\tau \leq \tau_\ast\).
  Moreover, \(U^L\) is locally smooth in $\xi$ and $\tau$.
\end{proposition}

The key difference from~\cite{kms} at this stage of the proof is that there, the
authors needed to solve the autonomous equation \(\partial_{\tau}V = \Tab V\),
for which it is trivial to construct a solution
\(U^L = \Re(e^{\tau\Tab}\eta) = \Re(e^{\tau\lambda}\eta)\) from an eigenvector
\(\eta \in E_{-\infty}\) of \(\Tab\). In the present situation, we find
ourselves with equation~\eqref{eq:linear-ancient}, which is non-autonomous due
to the $\tau$-dependence in the rescaled multiplier $M_\tau(|\nabla|)$, and
further the coefficient of the dilation operator $\xi \cdot \nabla_\xi$ also
depends on $\tau$. Firstly, we change variables again from $\tau$ to
$\theta = -2\alpha\log R(e^\tau)$ which has the effect of making the coefficient
of $\xi \cdot \nabla$ constant. After this, we need to construct the ancient
solution. Being a parabolic equation, we cannot simply invert and choose
\(U^L = \Sb( \theta_{0},\theta )^{- 1}\eta\). To get around this, we use the
invertibility of the semigroup only on the spectral subspace \(E_{-\infty}\) and
the convergence of the family to the semigroup as \(\theta \to - \infty\).

Fix \(\kappa \ll 1\) so that \(2\alpha - \kappa > 1\) and
\(2\alpha + \kappa < 3\). Define the change of time variable
\[
  \theta = \theta(\tau) \coloneqq - 2\alpha\log R( e^\tau ).
\]
Since
\begin{equation}
  \theta'(\tau) = 2\alpha b(\tau) \xrightarrow[\tau \to - \infty]{} 1\label{eq:theta-derivative}
\end{equation}
for \(\tau_\ast \ll - 1\), \(\theta( \cdot )\) is increasing on
\(\left( - \infty,\tau_\ast \right]\). Write \(\tau = \tau(\theta)\) for its
inverse, so that \(R( e^{\tau(\theta)} ) = e^{- \frac{\theta}{2\alpha}}\), and
set
\[
  c(\theta)
  \coloneqq \tau'(\theta)
  =         \frac{1}{2\alpha b( \tau(\theta) )} \xrightarrow[\theta\to-\infty]{} 1.
\]
In terms of \(\theta\), the equation~\eqref{eq:linear-ancient} is
\(\partial_{\theta}V = \Tbt V\), where
\begin{equation*}
  \begin{split}
    \Tzt & \coloneqq \left( c(\theta) - \frac{1}{2\alpha} \right)I + \frac{1}{2\alpha}\xi \cdot \nabla_\xi - c(\theta)\Mtt( |\nabla| ), \\
    \Tbt & \coloneqq \Tzt - c(\theta)\beta\Kop.
  \end{split}
\end{equation*}
The coefficient of \(\xi \cdot \nabla\) is now the constant
\(\frac{1}{2\alpha}\).

Let $(\mathbf A(\theta), D(\mathbf A(\theta)))_{\theta \in \Theta}$ be densely
defined closed operators that extend to operators
$\mathbf A(\theta): H^m(\mathbb R^3) \to \mathcal D'(\mathbb R^3)$. We recall
that \(( \Sop(\theta,r) )_{\theta \geq r}\) is an \emph{evolution family}
associated to \(\Aop(\theta)\) if it is a family of strongly continuous
operators with \(\mathbf S(r,r) = I\) and
\(\mathbf S(\theta,s)\mathbf S(s,r) = \mathbf S(\theta,r)\), and each
\(\mathbf S(\theta,r)u\) for \(u \in D( \Aop(r) )\) is a distributional solution
to \(\partial_{\theta}u = \Aop(\theta)u\). For the operator \(\Tbt\), the
evolution family $\Sb$ satisfies the Duhamel-type formula
\begin{equation}
  \Sb(\theta,r)
  = \Sz(\theta,r)
    - \beta\int_r^\theta\Sz(\theta,s)c(s)\Kop\Sb(s,r)ds.\label{eq:volterra-identity}
\end{equation}

\begin{lemma}
  \label{lem:fourier-multiplier-bound}

  Let \(\nu > 0\), \(d \geq 0\), and \(c > 0\). For \(0 < h \leq 1\),
  \begin{equation*}
    \sup\limits_{\rho \geq 0}( 1 + \rho^2 )^{\frac{d}{2}}e^{- ch\rho^\nu} \leq C_{c,d,\nu}h^{- \frac{d}{\nu}}.
  \end{equation*}
  Consequently, an operator of the form
  \[
    \widehat{T_{h}g}(\zeta) = \mu_{h}(\zeta)\hat{g}(a_h\zeta),
  \]
  where \(| {\mu_{h}(\zeta)} | \leq Ce^{- ch|\zeta|^\nu}\) and both \(a_h\) and
  \(a_h^{- 1}\) are uniformly bounded in $h$, maps \(H^s\) to \(H^{s+d}\) with
  \[
    \| T_h g\|_{H^{s+d}} \lesssim_s h^{-d/\nu} \|g\|_{H^s}.
  \]
\end{lemma}
\begin{proof}

  For \(0 \leq \rho \leq 1\), the left side is bounded by \(2^{\frac{d}{2}}\).
  For \(\rho \geq 1\), put \(y = (ch)^{\frac{1}{\nu}}\rho\). Then
  \[
    \rho^{d}e^{- ch\rho^\nu} =    (ch)^{- \frac{d}{\nu}}y^{d}e^{- y^\nu}
                             \leq C_{c,d,\nu}h^{- \frac{d}{\nu}},
  \]
  because the continuous function \(y^{d}e^{- y^\nu}\) is bounded on
  \(\left[ 0,\infty \right).\) Plancherel and a change of frequency variables
  give
  \[
    \begin{aligned}
      \|T_hg\|_{H^{s+d}}^2 & = \int_{\mathbb R^n}(1 + |\zeta|^2)^{s+d} |\mu_h(\zeta)|^2 |\widehat g(a_h\zeta)|^2\,d\zeta                                     \\
                           & \leq Ch^{-\frac{2d}{\nu}} \int_{\mathbb R^n}(1 + |\zeta|^2)^s |\widehat g(a_h\zeta)|^2\,d\zeta                                  \\
                           & = Ch^{-\frac{2d}{\nu}}|a_h|^{-n} \int_{\mathbb R^n}(1 + |a_h^{-1}\eta|^2)^s |\widehat g(\eta)|^2\,d\eta                         \\
                           & \leq C h^{-\frac{2d}{\nu}} \int_{\mathbb R^n}(1 + |\eta|^2)^s |\widehat g(\eta)|^2\,d\eta = C h^{-\frac{2d}{\nu}}\|g\|_{H^s}^2,
    \end{aligned}
  \]
  as claimed.
\end{proof}

\begin{lemma}
  \label{lem:volterra-evolution}

  For \(\theta_\ast \ll - 1\), \(\Tbt\) has an associated evolution family
  \(\Sb(\theta,r)\) on \(L_\sigma^2\) and on every \(H_\sigma^s\),
  \(s = 0,1,2,\ldots\) which agrees with the $L^2$ evolution. For every fixed
  \(H > 0\),
  \begin{equation}
    \| {\Sb(\theta,r)} \|_{H^s \to H^s} + (\theta - r)^{\frac{1}{2\alpha-\kappa}}\| {\Sb(\theta,r)} \|_{H^s \to H^{s+1}}
    \leq C_{s,H},\quad 0
    <    \theta - r
    \leq H.\label{eq:evolution-smoothing}
  \end{equation}
  In addition, for each fixed \(h > 0\),
  \begin{equation}
    \left\| {\Sb(\theta + h,\theta) - e^{h\Tab}} \right\|_{L^2 \to L^2}
    \to 0\quad\text{ as }\quad \theta
    \to - \infty.\label{eq:fixed-step-frozen-limit}
  \end{equation}
\end{lemma}

\begin{proof}

  First we understand the evolution for \(\beta = 0\). It turns out that the
  evolution family \(\Sz(\theta,r)\) can be given explicitly on the Fourier
  side:
  \begin{equation}
    ( \Sz(\theta,r)g )^\wedge(\zeta) = e^{\tau(\theta)-\tau(r)-2\frac{\theta-r}{\alpha}}e^{-(1-e^{\tau(r)-\tau(\theta)})\Mtt(|\zeta|)}\hat{g}(e^{\frac{-(\theta-r)}{2\alpha}}\zeta).\label{eq:free-propagator}
  \end{equation}
  Indeed, taking the Fourier transform of the equation
  \(\partial_{\theta}g = \Tzt g\), we find the transport equation
  \[
    \partial_\theta\hat{g}(\zeta) + \frac{1}{2\alpha}\zeta \cdot \nabla_\zeta\hat{g}(\zeta) = \left. c(\theta) - \frac{2}{\alpha} - c(\theta)\Mtt( |\zeta| ) \right.\hat{g}(\zeta).
  \]
  It follows from the explicit formula that \(\Sz\) is an evolution family.

  Fix \(H > 0\). Then \(1 - e^{\tau(r)-\tau(\theta)} \sim \theta - r\), while
  the Potter bounds in~\eqref{eq:potter-bounds} imply
  \(\Mtt(\rho) \geq c_0\rho^{2\alpha-\kappa} - C_0\). Thus
  Lemma~\ref{lem:fourier-multiplier-bound} gives
  \begin{equation}
    \| {\Sz(\theta,r)} \|_{H^s \to H^{s+j}} \leq C_{s,H}(\theta - r)^{- \frac{j}{2\alpha-\kappa}},\quad j
                                            =    0,1.\label{eq:free-smoothing}
  \end{equation}
  The semigroup \(e^{h\Taz}\) can be similarly shown to be given by
  \begin{equation}
    \left( e^{h\Taz}g \right)^\wedge(\zeta) = e^{\left( 1-\frac{2}{\alpha} \right)h}e^{-(1-e^{- h})|\zeta|^{2\alpha}}\hat{g}(e^{- \frac{h}{2\alpha}}\zeta).\label{eq:frozen-free-propagator}
  \end{equation}
  Hence, it is not hard to show using the bounds of
  Lemma~\ref{lem:symbol-calculus} that
  \begin{equation}
    \sup\limits_{h \in \lbrack h_{0},H\rbrack}\left\| {\Sz(\theta + h,\theta) - e^{h\Taz}} \right\|_{L^2 \to H^1} \xrightarrow[\theta \to - \infty]{} 0,\quad 0
    < h_0
    < H.\label{eq:free-frozen-limit}
  \end{equation}
  (Specifically one can use~\eqref{eq:regular-variation} on some
  \(\rho_0 \leq \rho \leq R_0\),~\eqref{eq:potter-bounds} on
  \(0 < \rho < \rho_0\), and bound by
  \(C(1 + \rho)e^{- c\rho^{2\alpha-\kappa}}\) for \(\rho > R_0\).)

  Now we use the Volterra series construction to obtain the evolution family for
  \(\Tbt\). Define the sequence \(( \Sin )_{n \geq 0}\) by
  \[
    \Siz
    \coloneqq \Sz,\quad\quad\Sio(\theta,r)
    \coloneqq \int_{r<r_{1}<\theta}\Sz( \theta,r_1 )\lbrack - c( r_1 )\beta\Kop \rbrack\Sz( r_{1},r )dr_{1},
  \]
  and for general \(n\),
  \begin{align}
    \Sin(\theta,r) & \coloneqq \int_{r<r'<\theta}\Sinm(\theta,r')\lbrack - c(r')\beta\Kop \rbrack\Sz(r',r)dr'                                                                                              \\
                   & = \int_{r<r'<\theta}\Sz(\theta,r')\lbrack - c(r')\beta\Kop \rbrack\Sinm(r',r)dr' \label{recursion}                                                                                    \\
                   & = \int_{r<r_{1}<\ldots<r_{n}<\theta}\Sz( \theta,r_n )\lbrack - c( r_n )\beta\Kop \rbrack\Sz( r_{n},r_{n-1} )\ldots\lbrack - c( r_1 )\beta\Kop \rbrack\Sz( r_{1},r )dr_1\ldots dr_{n}.
  \end{align}
  Each product \(\Kop\Sz\) is bounded by~\eqref{eq:free-smoothing}. Putting
  \(h \coloneqq \theta - r\), \(r_0 \coloneqq r\), and
  \(r_{n+1} \coloneqq \theta\), we find
  \begin{align}
    & \| {\Sin( r_{n+1},r_0 )} \|_{H^s \to H^{s+j}} \notag                                                                                                                                                                                                                                                  \\
    & \leq C^{n+1}\int_{r<r_{1}<\ldots<r_{n}<\theta}( r_{n+1} - r_n )^{- \frac{j}{2\alpha-\kappa}}\prod_{k=1}^{n}( r_k - r_{k-1} )^{- \frac{1}{2\alpha-\kappa}}dr_1\ldots dr_n\ \label{eq:volterra-simplex-majorant}\                                                                                       \\
    & = C^{n+1}\frac{\Gamma(1 - \frac{j}{2\alpha - \kappa}){\Gamma(1 - \frac{1}{2\alpha - \kappa})}^n}{\Gamma(n\left( 1 - \frac{1}{2\alpha - \kappa} \right) + 1 - \frac{j}{2\alpha - \kappa})}h^{n\left( 1-\frac{1}{2\alpha-\kappa} \right)-\frac{j}{2\alpha-\kappa}}.\ \label{eq:volterra-simplex-bound}\
  \end{align}
  The equality~\eqref{eq:volterra-simplex-bound} is because we recognize the
  Dirichlet simplex integral (put \(t_k = \frac{r_k - r_{k-1}}{r_{n+1} - r_0}\)
  so the ordered simplex is mapped to
  \(V_n = \left\{ t_k \geq 0:\sum_{k=1}^{n}t_k \leq 1 \right\}\), then
  consult~\cite[(5.14.2)]{NIST:DLMF}.) Here \(C\) may depend on \(s,H,\beta\)
  and the fixed background and symbol parameters, but not on \(n,\theta,r\). As
  $\Gamma(C_1n + C_2) \gg C^n\Gamma(C_1)^n$,
  \(\Sb \coloneqq \sum_{n \geq 0}\Sin\) converges locally in \(H^s \to H^s\)
  operator norm, and satisfies~\eqref{eq:evolution-smoothing}. Summing the
  recursive identity~\eqref{recursion} for \(\Sin\) gives the Duhamel
  formula~\eqref{eq:volterra-identity}.
  It follows from this that it defines solutions to
  \(\partial_{\theta}v = \Tbt v\).

  It remains to show~\eqref{eq:fixed-step-frozen-limit}. Let
  \(\Sab = \sum_{n \geq 0}\Sabn\) be the evolution family obtained by replacing
  \(\Sz(\theta,r)\) with \(e^{(\theta-r)\Taz}\) and \(c(s)\) with \(1\) in the
  above construction. We will show for fixed \(h > 0\) that
  \[
    \Sb(\theta + h,\theta) - \Sab(\theta + h,\theta) = \sum_{n \ge 0} (\Sin(\theta + h,\theta) - \Sabn(\theta + h,\theta) )\xrightarrow[\theta \to - \infty]{} 0
  \]
  as an operator on \(L^2\). Since clearly
  \(\Sab(\theta + h,\theta) = e^{h\Tab}\) by uniqueness, this gives us what we
  want. We only need to estimate the first \(M \gg 1\) terms as the tail sum is
  arbitrarily small by~\eqref{eq:volterra-simplex-bound} (which is satisfied by
  both $\Sin$ and $\Sabn$),
  \begin{equation}
    \label{eq:tail}
    \sum_{n \ge 0} (\Sin - \Sabn )
    = \sum_{n=0}^M (\Sin - \Sabn )
      + \underbrace{\sum_{n>M} (\Sin - \Sabn)}_{\xrightarrow[M\to\infty]{} 0}.
  \end{equation}
  For each fixed \(n \geq 1\), let
  \begin{align*}
    \Delta_n                  & \coloneqq \{(t_1,\ldots,t_n):0 = t_0 < t_1 < \ldots < t_n < t_{n+1} = h\},                                                                                         \\
    F_{\theta,n}(\mathbf{t})  & \coloneqq \Sz(\theta + h,\theta + t_n) [-c(\theta + t_n)\beta\Kop] \Sz(\theta + t_n,\theta + t_{n-1}) \cdots [-c(\theta + t_1)\beta\Kop] \Sz(\theta + t_1,\theta), \\
    F_{-\infty,n}(\mathbf{t}) & \coloneqq e^{(h-t_n)\Taz} [-\beta\Kop] e^{(t_n-t_{n-1})\Taz} \cdots [-\beta\Kop] e^{t_1\Taz},
  \end{align*}
  so that after putting $t_k = r_k - \theta$ we have
  \[
    \Sin(\theta + h,\theta)
    = \int_{\Delta_n} F_{\theta,n}(\mathbf{t})\,d\mathbf{t}, \qquad\text{and}\qquad \Sabn(\theta + h,\theta)
    = \int_{\Delta_n} F_{-\infty,n}(\mathbf{t})\,d\mathbf{t}.
  \]
  For every fixed \(\mathbf t \in \Delta_n\), all time increments $t_i - t_j$
  are strictly positive, so by~\eqref{eq:free-frozen-limit}, the boundedness of
  \(\mathbf K:H^1 \to L^2\), and \(c(\theta + t_k) \to 1\), we have
  \[
    \|F_{\theta,n}(\mathbf t)
    - F_{-\infty,n}(\mathbf t)\|_{L^2 \to L^2}\xrightarrow[\theta\to-\infty]{}0.
  \]
  Moreover, as in the proof of~\eqref{eq:volterra-simplex-majorant} for
  \(j = 0\),
  \[
    \|F_{\theta,n}(\mathbf t) - F_{-\infty,n}(\mathbf t)\|_{L^2 \to L^2} \leq 2C^{n+1}\prod_{k=1}^n(t_k - t_{k-1})^{-\frac1{2\alpha-\kappa}},
  \]
  which is integrable over \(\Delta_n\) by~\eqref{eq:volterra-simplex-bound}.
  Hence dominated convergence yields
  \[
    \|\mathbf S_{(n)}(\theta + h,\theta)
    - \mathbf S_{-\infty,(n)}(\theta + h,\theta)\|_{L^2 \to L^2}\xrightarrow[\theta\to-\infty]{}0.
  \]
  The case \(n = 0\) is just~\eqref{eq:free-frozen-limit}. Together
  with~\eqref{eq:tail}, this completes the proof.
\end{proof}

\begin{proof}[Proof of Proposition~\ref{prop:ancient-unstable-mode}]

  Let \(\theta_0 \ll - 1\), \(h_\ast \gg 1\),
  \(\theta_n = \theta_0 - nh_\ast \to - \infty\),
  \(e_\ast \in E_{-\infty} \setminus \{0\}\) real valued,
  \(\gamma_0 \in ( \sigma_{c},\sigma_u )\), and let us consider the problem of
  finding a sequence \(w = ( w_n )_{n \geq 0}\) of values for our solution
  \(w_n = V( \theta_n )\) solving
  \begin{equation}
    \begin{array}{rlr}
      w_n & = \Sb( \theta_{n},\theta_{n+1} )w_{n+1} \eqqcolon e^{h_\ast\Tab}w_{n+1} + \Dn w_{n+1},\quad\Pab w_0 & = e_\ast,
    \end{array}
    \label{eq:discrete-ancient-orbit0}
  \end{equation}
  and belonging to the weighted space \(\mathcal{X}_{\gamma_0}\) given by
  (recall \(\Pab\) is the projection \(L_\sigma^2 \to E_{-\infty}\))
  \begin{equation}
    \| w \|_{\mathcal{X}_{\gamma_0}}
    \coloneqq \sup\limits_{n \geq 0}e^{\gamma_{0}nh_\ast}\left( \| {\Pab w_n} \|_{L^2} + \| {( 1 - \Pab )w_n} \|_{L^2} \right)
    <         \infty.\label{eq:ancient-weighted-norm1}
  \end{equation}

  Given such a sequence, it is of course trivial to reconstruct a full solution
  \(V\) with the evolution family by setting
  \(V(\theta) = \Sb( \theta,\theta_{n+1} )w_{n+1}\) for
  \(\theta \in \lbrack \theta_{n+1},\theta_n \rbrack\). We will see later how
  the choices of the above parameters allow the construction to go through, and
  give \(V\) the required estimates. We can
  split~\eqref{eq:discrete-ancient-orbit0} into its projection to
  \(E_{-\infty}\) and the complement \(F_{-\infty}\):
  \begin{align*}
    \Pab w_{n+1}    & = e^{- h_\ast\Tab}\Pab w_n - e^{- h_\ast\Tab}\Pab\Dn w_{n+1},  \\
    ( I - \Pab )w_n & = e^{h_\ast\Tab}( I - \Pab )w_{n+1} + ( I - \Pab )\Dn w_{n+1}.
  \end{align*}
  where the operator \(e^{- h_\ast\Tab}\) makes sense because we are acting on
  the range of \(\Pab\). So let us define the backward and forward operators,
  \(\Bop \coloneqq \left( \left. e^{h_\ast\Tab}|_{E_{-\infty}} \right. \right)^{- 1}\)
  and \(\Fop \coloneqq e^{h_\ast\Tab}|_{F_{-\infty}}\). We iterate backwards to
  \(n = 0\) on \(E_{-\infty}\) and forwards to \(n \to \infty\) on
  \(F_{-\infty}\) to obtain the system
  \begin{align}
    \Pab w_n        & = \Bop^{n}e_\ast - \sum_{j=0}^{n-1}\Bop^{n-j}\Pab\Dj w_{j+1},\ \label{eq:half-line-green-a}\                                                     \\
    ( I - \Pab )w_n & = \sum_{j=n}^\infty\Fop^{j-n}( I - \Pab )\Dj w_{j+1} + \lim\limits_{J \to \infty}\Fop^{J-n+1}( I - \Pab )w_{J+1}.\ \label{eq:half-line-green-b}\
  \end{align}
  By Proposition~\ref{prop:frozen-dichotomy},
  \(\gamma_0 \in ( \sigma_{c},\sigma_u )\), and \(h_\ast \gg 1\), we have that
  \[
    \max(e^{\gamma_{0}h_\ast}\| \Bop \|,e^{- \gamma_{0}h_\ast}\| \Fop \|)
    \eqqcolon a
    <         1.
  \]
  In particular, this second term of~\eqref{eq:half-line-green-b} is in fact
  zero:
  \[
    \| {\Fop^{J-n+1}( I - \Pab )w_{J+1}} \|_{L^2} \leq a^{J-n+1}e^{- \gamma_{0}nh_\ast}\| w \|_{\mathcal{X}_{\gamma_0}} \xrightarrow[J \to \infty]{} 0.
  \]
  By Lemma~\ref{lem:volterra-evolution}, \(\theta_0 \ll - 1\) gives us a
  \(\delta\) so that
  \begin{equation}
    e^{-\gamma_0h_\ast}\sup_{n \ge 0} \|\Dn\|_{L^2 \to L^2}
    \le \delta, \qquad \delta\, \frac{a\|\Pab\| + \|I - \Pab\|}{1 - a}
    <   1. \label{eq:weighted-small-perturbation0}
  \end{equation}
  It follows that we can take the right-hand sides
  of~\eqref{eq:half-line-green-a} and~\eqref{eq:half-line-green-b} as defining
  an affine contraction \(\mathcal{X}_{\gamma_0} \to \mathcal{X}_{\gamma_0}\),
  giving the existence of a fixed point \(w\), and hence the existence of \(V\)
  satisfying the estimates
  \(\| {V(\theta)} \|_{L^2} \leq Ce^{\gamma_{0}(\theta-\theta_{0})}.\) By
  using~\eqref{eq:evolution-smoothing} repeatedly we can upgrade this to
  \[
    \| {V(\theta)} \|_{H^k} \leq C_k\| {V(\theta - 1)} \|_{L^2}
                            \leq C_{k}e^{\gamma_{0}(\theta-\theta_{0})}.
  \]
  Since $\Pab V(\theta_0) = e_\ast \ne 0$, we have $V(\theta_0) \ne 0$. The
  construction preserves real-valued fields. Moreover, $V(\theta)$ cannot vanish
  at an earlier time, since the evolution property would then imply
  $V(\theta_0) = 0$. To finish the proof we undo the change of variables
  \(\tau \to \theta\), choose $\tau_0$ so that $\theta_0 = \theta(\tau_0)$, and
  rescale,
  \[
    U^{L}(\tau) \coloneqq \frac{V( \theta(\tau) )}{\| {V( \theta_0 )} \|_{L^2}}.
  \]
  By taking \(\theta_0 \ll - 1\) and using~\eqref{eq:theta-derivative}, we get
  that \(( \theta - \theta_0 ) \leq (1 - \varepsilon)( \tau - \tau_0 )\) for
  some \(\varepsilon \ll 1\) which we may choose so that
  \(\gamma \coloneqq \gamma_{0}(1 - \varepsilon) > \frac{1}{2}\). Then the
  estimates for \(V\) implies~\eqref{eq:UL-bounds} as needed.

  Finally, we need to prove $U^L$ is locally smooth in spacetime. We do this by
  converting back to physical coordinates $(x,t)$ to avoid estimating terms with
  the operator $\xi \cdot \nabla$. Set
  \[
    u^L(x,t) = \frac1{t R(t)} U^L(R(t)x, \log t), \quad \overline u(x,t)
             = \frac1{t R(t)} \overline U(R(t)x).
  \]
  The equation for $u^L$ is
  \[
    \partial_t u^L
    = - m(|\nabla|) u^L
      - \beta \PP(\overline u \cdot \nabla u^L + u^L \cdot \nabla \overline u).
  \]
  Since $u^L \in C^0 H^k$ locally in $(t,x)$, this equation can be used to
  inductively show that $u^L \in C^j H^k$ for $j = 1,2,3,\dots$, which finishes
  the proof after changing back to $(\xi,\tau)$ variables.
\end{proof}

\section{Leray--Hopf nonuniqueness}
\label{sec:construction}

\begin{proof}[Proof of Theorem~\ref{thm:main}]

  Put
  \begin{equation}
    F_\beta \coloneqq  {\PP} \lbrack \beta(b(\tau) - 1 - b(\tau)\xi \cdot \nabla_\xi)\overline{U} + \beta \Mt( |\nabla| )\overline{U} + \beta^2\overline{U} \cdot \nabla_\xi\overline{U} \rbrack,\quad U^{( \pm )}
            \coloneqq        \beta\overline{U}
                             \pm U^{L}.\label{eq:golovkin-profiles}
  \end{equation}

  This is Golovkin's trick~\cite{golovkin-1964, golovkin-1965}, as explained by
  Dolce--Mescolini~\cite[\S 1.2]{dolce-mescolini}: the terms linear in \(U^L\)
  cancel by~\eqref{eq:linear-ancient}, while
  \({\PP}(( \pm U^L ) \cdot \nabla_\xi( \pm U^L )) = {\PP}(U^L \cdot \nabla_{\xi}U^{L})\).
  Hence, both \(U^\pm\) solve the projected equation with the same force
  \(F_\beta + {\PP}(U^L \cdot \nabla_{\xi}U^{L})\).

  Set \(t_\ast \coloneqq e^{\tau_\ast}\) and define
  \begin{align}
    u^{( \pm )}(x,t) & \coloneqq \frac{1}{tR(t)}U^{( \pm )}( R(t)x,\log t ),                                                                                       \\
    f(x,t)           & \coloneqq \frac{1}{t^{2}R(t)}\lbrack F_\beta + {\PP}(U^L \cdot \nabla_{\xi}U^{L}) \rbrack( R(t)x,\log t ).\label{eq:physical-golovkin-pair}
  \end{align}
  These fields are real, distinct, and divergence-free. Moreover, they are
  smooth for \(t > 0\), belonging to every \(H^k\) uniformly on compact subsets
  of \((0,t_\ast]\). Hence, \(u = u^{( \pm )}\) solve the
  equation~\eqref{eq:main}. It remains to verify that
  \(u^{( \pm )} \in L^{\infty}L^2 \cap L^2{\dot{H}}_m\), and that
  \(f \in L^2\mathcal{V}'_m\), with \(f \in L^{1}L^2\) if
  \(\alpha < \frac{5}{4}\) or \(\alpha = \frac{5}{4}\) and
  \(\ell \in L^1\left( (1,\infty);\frac{dr}{r} \right)\), and the energy
  equality~\eqref{eq:forced-energy} holds at the initial time.

  The properties of \(\overline{U}\) recalled in
  Section~\ref{sec:ancient-instability} and~\eqref{eq:UL-bounds} give
  \(\sup_{\tau \leq \tau_\ast}\| {U^{( \pm )}(\tau)} \|_{H^2} \leq C,\) and
  hence
  \begin{equation}
    \| {u^{( \pm )}(t)} \|_{L_x^2}
    =    \frac{\| {U^{( \pm )}( \log t )} \|_{L_\xi^2}}{t{R(t)}^{\frac{5}{2}}}
    \leq C{R(t)}^{- \frac{5}{2}}m( R(t) ) \xrightarrow[t \downarrow 0]{} 0.\label{eq:construction-trace}
  \end{equation}
  So \(u^{( \pm )} \in L^{\infty}L^2\) and the zero initial data is attained. To
  see that \(u^{( \pm )} \in L^2{\dot{H}}_m\), first observe the inequality
  \[
    \int_{\mathbb{R}^3}\Mt( |\eta| )| {\hat{W}(\eta)} |^{2}d\eta \leq C\| W \|_{H^2}^{2}.
  \]
  which holds as \(\Mt( |\eta| ) \lesssim \begin{cases}
                                            |\eta|^{2\alpha-\kappa}, & 0 < |\eta| \le 1, \\
                                            |\eta|^{2\alpha+\kappa}, & |\eta| \ge 1,
                                          \end{cases}\) for $\kappa > 0$, and
  \(\alpha \leq \frac{5}{4} < 2\) for \(\tau \ll 1\). Thus,
  \begin{align}
    \| {u^{( \pm )}(t)} \|_{{\dot{H}}_m}^2 & = \frac{1}{t^2{R(t)}^5}\int_{\mathbb{R}^3}m( R(t)|\eta| )| {{\hat{U}}^{( \pm )}(\eta,\tau)} |^{2}d\eta \\
                                           & = \frac{1}{t^3{R(t)}^5}\int_{\mathbb{R}^3}\Mt( |\eta| )| {{\hat{U}}^{( \pm )}(\eta,\tau)} |^{2}d\eta   \\
                                           & \leq \frac{C}{t}({R(t)}^{- \frac{5}{2}}m( R(t) ))^2 . \label{eq:ineq-integrable}
  \end{align}
  Then we integrate in \(t\) and change variables \(r = R(t)\). Since
  differentiating \(m( R(t) ) = t^{- 1}\) gives
  \(\frac{dt}{t} = - \frac{Rm'(R)}{m(R)}\frac{dR}{R}\),
  \(\frac{Rm'(R)}{m(R)} \to 2\alpha,\) we obtain
  \begin{equation}
    \int_0^{t_\ast}\| {u(t)} \|_{{\dot{H}}_m}^{2}dt
    \leq C\int_{R(t_\ast)}^\infty\left( r^{- \frac{5}{2}}m(r) \right)^2\frac{dr}{r}
    <    \infty.\label{eq:construction-dissipation}
  \end{equation}

  Next we estimate \(f\). For \(\frac{1}{t^{2}R(t)}F_\beta( R(t)x,\log t ),\)
  splitting the below integral into \(|\eta| \lesssim 1\) (where
  \(\Mt(\eta) \gtrsim |\eta|^{2\alpha+\kappa}\)) and \(|\eta| \gtrsim 1\) (where
  \(\Mt(\eta) \gtrsim 1\)) gives
  \[
    \int_{\mathbb{R}^3}\frac{| {{\hat{F}}_\beta(\eta,\tau)} |^2}{1 + m( R(t)|\eta| )}d\eta \leq \frac{C}{m( R(t) )}.
  \]

  Hence,
  \begin{align*}
    \left\| {\frac{1}{t^{2}R(t)}F_\beta( R(t)x,\log t )} \right\|_{\mathcal{V}_{m}'}^2 & \leq \frac{1}{t^4{R(t)}^5}\int_{\mathbb{R}^3}\frac{| {{\hat{F}}_\beta( \eta,\log t )} |^2}{1 + m( R(t)|\eta| )}d\eta \\
                                                                                       & \leq \frac{C}{t^4{R(t)}^{5}m( R(t) )} = \frac{C}{t}\left( {R(t)}^{- \frac{5}{2}}m( R(t) ) \right)^{2},
  \end{align*}
  which is the same upper bound in~\eqref{eq:ineq-integrable}, so this is
  integrable in time.

  For the remaining part
  \(\frac{1}{t^{2}R(t)}{\PP}(U^L \cdot \nabla_{\xi}U^{L})( R(t)x,\log t ),\)
  \[
    \left\| {\frac{1}{t^{2}R(t)}{\PP}(U^L \cdot \nabla_{\xi}U^{L})( R(t)x,\log t )} \right\|_{L^2}^2
    \leq C{R(t)}^{- 5}t^{4\gamma-4}
    =    Ct^{4\gamma-2}\left( R^{- \frac{5}{2}}m(R) \right)^{2}.
  \]
  Since \(R^{- \frac{5}{2}}m(R)\) is bounded by~\eqref{eq:symbol-little-o} and
  \(\gamma > \frac{1}{2}\), this term belongs to
  \(L^{2}( ( 0,t_\ast );L^2 ) \subset L^{2}( ( 0,t_\ast );\mathcal{V}'_m ).\)
  This proves \(f \in L^2\mathcal{V}'_m\) as needed. The estimates proven above
  allow us to send \(s \downarrow 0\) in the energy equality.

  The last thing we need to show is that under \(\alpha < \frac{5}{4}\) or
  \(\alpha = \frac{5}{4}\) and
  \(\ell \in L^1\left( (1,\infty);\frac{dr}{r} \right)\), \(f \in L^{1}L^2\)
  also. The quadratic part is \(L^{2}L^2\) so we only need to deal with the part
  from \(F_\beta\). As \(\| F_\beta \|_{L^2} \leq C\), we obtain
  \[
    \left\| {\frac{1}{t^{2}R(t)}F_\beta( R(t)x,\log t )} \right\|_{L^2}
    \leq \frac{C}{t^2{R(t)}^{\frac{5}{2}}}
    =    \frac{C}{t}{R(t)}^{- \frac{5}{2}}m( R(t) ).
  \]
  Integrating and changing variables, we find that
  \[
    \int_0^{t_\ast}\| {f_\beta(t)} \|_{L^2}dt \leq C\int_{R(t_\ast)}^{\infty}r^{- \frac{5}{2}}m(r)\frac{dr}{r},
  \]
  which is finite under our assumptions as needed. This finishes the proof.
\end{proof}



\section{Uniqueness}
\label{sec:uniqueness}

In the \(\alpha > \frac{5}{4}\) case, the classical argument closes already with
Gronwall's inequality at the \(L^2\) level, because
\(1 + m(r) \gtrsim \langle r \rangle^{\frac{5}{2}}\) i.e.
\(\mathcal{V}_m \hookrightarrow H^{\frac{5}{4}}\), so we have the estimate
\[
  |\langle u \cdot \nabla v,\, w \rangle|
  \leq     \| u \|_{L^2}\| {\nabla v} \|_{L^{\frac{12}{5}}}\| w \|_{L^{12}}
  \lesssim \| u \|_{L^2}\| {v} \|_{\mathcal{V}_m}\| w \|_{\mathcal{V}_m}.
\]
We omit the details, which can for example be found in~\cite{hussein}
\emph{mutatis mutandis}. Henceforth we set \(\alpha = \frac{5}{4}\). In this
case we may not have \(\mathcal{V}_m \hookrightarrow H^{\frac{5}{4}}\), so we
repair this by studying the uniqueness in the negative regularity space
\(H^{- 1}\).

\begin{lemma}
  \label{lem:uniqueness-nonlinear-est}

  For \(u,v \in \mathcal{V}_m\),
  \begin{equation}
    \| {\langle \nabla \rangle\ell(\langle \nabla \rangle)(u \otimes v)} \|_{L^2} \lesssim \| u \|_{\mathcal{V}_m}\| v \|_{\mathcal{V}_m}.\label{eq:uniqueness-bd-1}
  \end{equation}
  Consequently, for vector fields \(u \in \mathcal{V}_m\) and \(v \in L^2\),
  \begin{equation}
    \| {\langle \nabla \rangle^{- 1}{\PP}\nabla \cdot (u \otimes v)} \|_{\mathcal{V}'_m} \lesssim \| u \|_{\mathcal{V}_m}\| {{\ell(\langle \nabla \rangle)}^{- 1}\langle \nabla \rangle^{- 1}v} \|_{L^2}.\label{eq:uniqueness-bd-2}
  \end{equation}
\end{lemma}
\begin{proof}
  Let \(K(\zeta) = \frac{1}{1 + m( |\zeta| )}\). By Cauchy--Schwarz, we have
  \begin{equation}
    | {\widehat{u_{i}v_{j}}(\zeta)} |^2 \lesssim (K \ast K)(\zeta)\int_{{\mathbb{R}}^3}( 1 + m( |\eta| ) )( 1 + m( |\zeta - \eta| ) )|\hat{u}_{i}(\eta)|^{2}|\hat{v}_{j}(\zeta - \eta)|^{2}d\eta.\label{ineq:KstarK}
  \end{equation}
  Hence, to prove~\eqref{eq:uniqueness-bd-1}, it suffices by Plancherel and
  Fubini to show that
  \begin{equation}
    K \ast K(\zeta) \lesssim \langle \zeta \rangle^{- 2}{\ell(\langle \zeta \rangle)}^{- 2}.\label{ineq:uniqueness-conv}
  \end{equation}
  For \(|\zeta| \leq 4\) , \(|K \ast K(\zeta)| \lesssim 1\) by Holder's
  inequality, so we may assume \(|\zeta| > 4.\) Karamata's integral
  theorem~\cite[Theorem 1.5.11(i)–(ii)]{bingham-goldie-teugels} gives for
  \(R \geq 2\) that
  \begin{equation}
    \int_{|\eta| \leq R}K(\eta)d\eta
    \lesssim \frac{R^{\frac{1}{2}}}{\ell(R)},\quad\int_{|\eta| \geq R}{K(\eta)}^{2}d\eta
    \lesssim \frac{R^{- 2}}{{\ell(R)}^2},\label{ineq:integral-K-bounds}
  \end{equation}
  as the integrands in radial coordinates are \(\frac{r^2}{1 + m(r)}\) and
  \(\frac{r^2}{( 1 + m(r) )^2}\), which vary regularly with indices
  \(- \frac{1}{2} > - 1\) and \(- 3 < - 1\) respectively.

  We split
  \(\int_{{\mathbb{R}}^3}K(\eta)K(\zeta - \eta)d\eta = I_1 + I_2 + I_3\), where

  \begin{gather*}
    I_i = \int_{U_i}K(\eta)K(\zeta - \eta)d\eta,\quad\quad U_1 = \left\{ \eta \in {\mathbb{R}}^3 :|\eta| \leq \frac{|\zeta|}{2} \right\},\quad\quad U_2 = \left\{ \eta \in {\mathbb{R}}^3 :|\zeta - \eta| \leq \frac{|\zeta|}{2} \right\}, \\
    U_3 = {\mathbb{R}}^3\backslash( U_1 \cup U_2 ) = \left\{ \eta \in {\mathbb{R}}^3 :|\eta| > \frac{|\zeta|}{2}\text{ and }|\zeta - \eta| > \frac{|\zeta|}{2} \right\}.
  \end{gather*}
  For \(I_1\), we have
  \(\frac{|\zeta|}{2} \leq | {\zeta - \eta} | \leq \frac{3}{2}|\zeta|\),
  so~\eqref{eq:potter-bounds} gives that
  \(|K(\zeta - \eta)| \lesssim \frac{1}{|\zeta|^{\frac{5}{2}}\ell(|\zeta|)}\).
  Hence, by~\eqref{ineq:integral-K-bounds}, we have
  \(|I_{1}| \lesssim |\zeta|^{- 2}{\ell(|\zeta|)}^{- 2}.\) \(I_2\) is treated in
  the same way by symmetry. For the remaining \(I_3\) we have
  by~\eqref{ineq:integral-K-bounds} and Cauchy-Schwarz that
  \[
    |I_{3}|
    =        \left| {\int_{|\eta|>\frac{|\zeta|}{2}\text{ and }|\zeta-\eta|>\frac{|\zeta|}{2}}K(\eta)K(\zeta - \eta)d\eta} \right|
    \leq     \int_{|\eta|>\frac{|\zeta|}{2}}{K(\eta)}^{2}d\eta
    \lesssim \frac{|\zeta|^{- 2}}{{\ell(|\zeta|)}^2}.
  \]
  This shows~\eqref{eq:uniqueness-bd-1}. The
  inequality~\eqref{eq:uniqueness-bd-2} follows from this, as if \(\varphi\) is
  a test function,
  \begin{align*}
    \left| {\int_{{\mathbb{R}}^3}\langle \nabla \rangle^{- 1}{\PP}\nabla \cdot (u \otimes v)\varphi dx} \right| & = \left| {\int_{{\mathbb{R}}^3}u \otimes v\nabla{\PP}\langle \nabla \rangle^{- 1}\varphi dx} \right|                                                                                                                 \\
                                                                                                                & \lesssim \| {{\ell(\langle \nabla \rangle)}^{- 1}\langle \nabla \rangle^{- 1}v} \|_{L^2}\| {\ell(\langle \nabla \rangle)\langle \nabla \rangle( u \otimes \nabla{\PP}\langle \nabla \rangle^{- 1}\varphi )} \|_{L^2} \\
                                                                                                                & \lesssim \| {{\ell(\langle \nabla \rangle)}^{- 1}\langle \nabla \rangle^{- 1}v} \|_{L^2}\| u \|_{\mathcal{V}_m}\| \varphi \|_{\mathcal{V}_m}.\qedhere
  \end{align*}
\end{proof}

\begin{proof}[Proof of Theorem~\ref{thm:local-wellposedness}]

  Let \(u^{(1)},u^{(2)}\) be two solutions on the time interval
  \(\lbrack 0,T\rbrack\) with the same initial data and forcing, and let
  \(z = \langle \nabla \rangle^{- 1}( u^{(1)} - u^{(2)} ) \in L^{2}( 0,T;\mathcal{V}_m )\).
  It solves the equation
  \begin{equation}
    \partial_{t}z + m( |\nabla| )z + \langle \nabla \rangle^{- 1}{\PP}\nabla \cdot ( u^{(1)} \otimes \langle \nabla \rangle z + \langle \nabla \rangle z \otimes u^{(2)} ) = 0.\label{eq:z}
  \end{equation}
  Observe that
  \(\| {m( |\nabla| )z} \|_{\mathcal{V}'_m} \lesssim \| z \|_{{\dot{H}}_m}\),
  and for the nonlinear term we have by~\eqref{eq:uniqueness-bd-2} that
  \[
    \| {\langle \nabla \rangle^{- 1}{\PP}\nabla \cdot ( u_i^{(k)}\langle \nabla \rangle z_j )} \|_{L^2\mathcal{V}'_m}
    \lesssim \| u_i^{(k)} \|_{L^2\mathcal{V}_m}\| {{\ell(\langle \nabla \rangle)}^{- 1}z_j} \|_{L^{\infty}L^2}
    \lesssim \| u_i^{(k)} \|_{L^2\mathcal{V}_m}\| {u^{(1)} - u^{(2)}} \|_{L^{\infty}L^2}.
  \]
  Hence, we see that
  \[
    \partial_{t}z \in L^{2}( 0,T;\mathcal{V}'_m ).
  \]
  Thus, \(z \in C^{0}L_\sigma^2\) with \(z(0) = 0\) (see
  e.g.~\cite[Chapter~III]{temam}). Testing~\eqref{eq:z} against \(z\) and
  using~\eqref{eq:uniqueness-bd-2} again with
  \(\| z \|_{\mathcal{V}_m}^2 \simeq \| z \|_{L^2}^2 + \left\| {{m( |\nabla| )}^{\frac{1}{2}}z} \right\|_{L^2}^2\),
  and applying Young's inequality for products, we have
  \[
    \begin{array}{r}
      \partial_t\| z \|_{L^2}^2 + \left\| {{m( |\nabla| )}^{\frac{1}{2}}z} \right\|_{L^2}^2 \leq C\left( \| u^{(1)} \|_{\mathcal{V}_m} + \| u^{(2)} \|_{\mathcal{V}_m} \right)\| {{\ell(\langle \nabla \rangle)}^{- 1}z} \|_{L^2}\| z \|_{\mathcal{V}_m} \\
                   \leq \frac{\| {{m( |\nabla| )}^{\frac{1}{2}}z} \|_{L^2}^2}{2} + \frac{\| z \|_{L^2}^2}{2} + C\left( \| u^{(1)} \|_{\mathcal{V}_m}^2 + \| u^{(2)} \|_{\mathcal{V}_m}^2 \right)\| {{\ell(\langle \nabla \rangle)}^{- 1}z} \|_{L^2}^{2}.
    \end{array}
  \]
  This first term on the right can be absorbed into the left hand side. So to
  close the estimate we need to estimate
  \(\| {{\ell(\langle \nabla \rangle)}^{- 1}z_j} \|_{L^2}^2\). Let \(R > 1\) be
  a parameter to be chosen shortly. We write
  \[
    \| {{\ell(\langle \nabla \rangle)}^{- 1}z} \|_{L^2}^2 = \left( \int_{\langle\zeta\rangle<R} + \int_{\langle\zeta\rangle \geq R} \right){\ell(\langle \zeta \rangle)}^{- 2}| {\hat{z}(\zeta)} |^{2}d\zeta.
  \]
  For \(\langle \zeta \rangle \leq R\), we have
  \({\ell(\langle \zeta \rangle)}^{- 2} \leq {\underline{\ell}(R)}^{- 2}\) by
  the definition of \(\underline{\ell}\). For \(\langle \zeta \rangle > R\), we
  estimate using~\eqref{eq:potter-bounds},
  \[
    {\ell(\langle \zeta \rangle)}^{- 2}
    \leq {\ell(R)}^{- 2}             \frac{\langle \zeta \rangle^2}{R^2}
    \leq {\underline{\ell}(R)}^{- 2} \frac{\langle \zeta \rangle^2}{R^2}.
  \]
  Hence, we find that
  \[
    \| {{\ell(\langle \nabla \rangle)}^{- 1}z} \|_{L^2}^2 \leq \frac{\| z \|_{L^2}^2 + R^{- 2}\| {\langle \nabla \rangle z} \|_{L^2}^2}{{\underline{\ell}(R)}^2}
  \]
  Let
  \(B = 1 + \| {u^{(1)} - u^{(2)}} \|_{L^\infty(0,T;L^{2})}^2 \geq \| {\langle \nabla \rangle z(t)} \|_{L^2}^2\),
  and choose \(R^2 = \frac{B}{\| z \|_{L^2}^2}\) to find the inequality
  \[
    \partial_t\| {z(t)} \|_{L^2}^2 \leq C\left( 1 + \| {u^{(1)}(t)} \|_{\mathcal{V}_m}^2 + \| {u^{(2)}(t)} \|_{\mathcal{V}_m}^2 \right)\frac{\| {z(t)} \|_{L^2}^2}{\underline{\ell}(\frac{\sqrt{B}}{\| {z(t)} \|_{L^2}})^2}.
  \]
  Since \(\| {u^{(i)}(\cdot)} \|^2_{\mathcal{V}_m} \in L^1([0,T])\), and
  $Z \mapsto \frac Z{\underline{\ell}(\sqrt{B/Z})^2}$ is an Osgood modulus, i.e.
  for any \(0 < c < B\),
  \[
    \int_0^c{\underline{\ell}\bigg( \sqrt{\frac{B}{Z}} \bigg)}^2\frac{dZ}{Z}
    = 2\int_{\sqrt{\frac{B}{c}}}^\infty{\underline{\ell}(r)}^2\frac{dr}{r}
    = \infty,
  \]
  it follows from the Osgood lemma~\cite[Lemma 3.4]{bcd} that \(z = 0\), and
  hence \(u^{(1)} = u^{(2)}\) as claimed.
\end{proof}

\section{
  Declaration of generative AI and AI-assisted technologies in the manuscript
  preparation process
}

This work was created with the help of OpenAI's ChatGPT5.5 and 5.6 which helped
clarify ideas. The paper itself was completely written by hand, followed by
further checks with ChatGPT for mathematical, typographical, and bibliographical
errors which were verified manually before being used to improve the paper. We take full responsibility for its contents.
\printbibliography

\end{document}